\documentclass[12pt,reqno]{amsart}
\usepackage{color,xcolor,pifont,mathrsfs,tikz}
\usepackage{hyperref}
\hypersetup{
    colorlinks=true,       
    linkcolor=blue,          
    citecolor=blue,        
    filecolor=blue,      
    urlcolor=cyan  
}
\usepackage{latexsym,amssymb,amsmath,amsthm}
\newcommand{\R}{{\ensuremath{\mathbb{R}}}}
\newcommand{\N}{{\ensuremath{\mathbb{N}}}}
\newcommand{\Z}{{\ensuremath{\mathbb{Z}}}}

\renewcommand{\P}{\ensuremath{\mathbb{P}}}
\renewcommand{\dj}{d\kern-0.4em\char"16\kern-0.1em}
\newcommand{\E}{\ensuremath{\mathbb{E}}}

\newcommand{\sig}{\ensuremath{\mathcal{F}}}

\newcommand{\wt}[1]{{\widetilde{#1}}}

\newtheorem{Thm}{Theorem}[section]
\newtheorem{Cor}[Thm]{Corollary}
\newtheorem{Lem}[Thm]{Lemma}
\newtheorem{Prop}[Thm]{Proposition}

\theoremstyle{remark}
\newtheorem{Rem}[Thm]{Remark}

\theoremstyle{definition}
\newtheorem{Ex}[Thm]{Example}

\theoremstyle{definition}

\theoremstyle{definition}

\begin{document}
\numberwithin{equation}{section}
\bibliographystyle{amsalpha}

\title[Heat kernel estimates for SBM]{Heat kernel estimates for subordinate Brownian motions}
\begin{abstract}
In this article we study transition probabilities of a class of subordinate Brownian motions. Under mild assumptions on the Laplace exponent of the corresponding subordinator, sharp two sided estimates of the transition probability are established. This approach, in particular, covers subordinators with Laplace exponents that vary regularly at infinity with index one, e.g. 
\[
 \phi(\lambda)=\frac{\lambda}{\log(1+\lambda)}-1 \quad \text{ or }\quad \phi(\lambda)=\frac{\lambda}{\log(1+\lambda^{\beta/2})},\ \beta\in (0,2)\,
\]
that correspond to subordinate Brownian motions with scaling order that is not necessarily stricty between $0$ and $2$. These estimates are applied to estimate Green function (potential) of  subordinate Brownian motion. We also prove the equivalence of the lower scaling condition of the Laplace exponent and the near diagonal upper estimate of the transition estimate. 
\end{abstract}

 \author{Ante Mimica}
 \address{Department of Mathematics, University of Zagreb, Bijeni\v cka cesta 30, 10 000 Zagreb, Croatia}
 \curraddr{}
 \thanks{Supported in part by Croatian Science Foundation under the project 3526.}
 \email{amimica@math.hr}



\maketitle

\allowdisplaybreaks[3]

\section{Introduction}

Recently, in probability theory and analysis there has been made progress in study of various properties of discontinuous Markov processes and their associated non-local generators. One of the notions that connects these two subjects is the heat kernel. In probability theory it can be understood as the transition density $p(t,x,y)$ of a Markov process $X$, while in analysis it is the fundamental solution of the equation $\partial_t u =\mathcal{A}u$, where $\mathcal{A}$ is the infinitesimal generator of $X$. Hence it is not surprising that one of the problems that drew much attention recently was to find sharp estimates of the heat kernel $p(t,x,y)$ for various classes of discontinuous Markov processes $X$ and non-local operators $\mathcal{A}$. 

For pure jump symmetric processes with stable-like jumping kernels on $\Z^d$ or $\R^d$ sharp heat kernel estimates were obtained in \cite{BL2,CK2,CK}, while in \cite{BGR,Kn, KSz2, KSz,Sz3} they were obtained for some classes of L\' evy processes. 
A common property of all stochastic processes for which sharp two-sided estimates of the heat kernel were obtained is that the scaling order was always strictly between $0$ and $2$. 
This notion will be explained in detail later (see (\ref{eq:lowerscaling}) and (\ref{eq:upperscaling}); for alternative approach to scaling Matuszewska indices may also be used, see \cite{BGT}).  Our motivation was to obtain sharp heat kernel estimates of heat kernel when this property fails. 
Some far from optimal upper bounds in this case were obtained earlier in \cite{Mi3} and \cite{KSz}.

Markov jump processes became important also in applications (e.g. in physics and  finance, see \cite{CT1}) and typical examples are $\alpha$-stable processes, where $\alpha\in (0,2)$ and from this aspect it is also important to have good estimates of transition densities. Processes that were not covered by the theory known so far are conjugate geometric stable processes. These are L\' evy processes $X=(X_t)_{t\geq 0}$ in $\R^d$ such that, for some $\beta\in (0,2]$,
\[
 \E e^{i\xi\cdot X_t}=e^{-t\frac{|\xi|^2}{\log(1+|\xi|^\beta)}},\qquad \xi\in \R^d\,.
\]
The order of such processes is not strictly less than $2$. Actually,  concerning behavior of jumps and some other potential-theoretic notions (e.g. Green function or jumping kernel), these processes are between any rotationally invariant $\alpha$-stable process ($\alpha\in (0,2)$) and Brownian motion. It can be seen (see (\ref{eq:green-conjgamma}) and \cite{Mi2,Mi3}) that the intensity of small jumps is higher than in the case of any stable process.  Knowing sharp estimates of transition density of such processes might be useful in modeling various phenomena in nature by them. 

A very successful technique that was used to obtain upper bounds in heat kernel estimates was developed by Carlen, Kusuoka and Stroock in the paper \cite{CKS} and some of the already mentioned papers actually use this method. This method works well in stable-like cases (see  \cite{BGK,BL2,CK,CKK,CKK2,KSch}), but it is not clear how to extend it to the cases when the scaling order is not strictly between $0$ and $2$. Hence, obtaining sharp heat kernel estimates for subordinate Brownian motion not satisfying this scaling order restriction could be a possible starting point for developing a generalized version of this method. 

Another motivation for this investigation was to try to generalize existing heat kernel estimates within the class of subordinate Brownian motions. As it will be seen from the main result, estimates will have a new form, but when the scaling order is restricted to be strictly between $0$ and $2$, these estimates reduce to the already known form.  

For $a,b\in \R$ we denote $a\wedge b:=\min\{a,b\}$ and $a\vee b:=\max\{a,b\}$. Notation $f(x)\asymp g(x), x\in I$ means that there exist constants $c_1,c_2>0$ such that $c_1f(x)\leq g(x)\leq c_2g(x)$ for $x\in I$. By $B_r(x)=\{y\in\R^d: |x-y|<r\}$ we denote open ball around $x\in \R^d$ with radius $r>0$\,. We also use convention $0^{-1}=+\infty$. 

A subordinator $S=(S_t)_{t\geq 0}$ is an increasing L\' evy process, which is a stochastic process defined on a probability space $(\Omega,\sig,\P)$ with stationary and independent increments with sample paths that are right continuous with left limits. It follows from the definition that $S$ takes values in $[0,\infty)$ and the Laplace transform of $S_t$ is of the form
\[
    \E e^{-\lambda S_t}=e^{-t\phi(\lambda)},\qquad \lambda>0,
\]
where $\phi$ is called the Laplace exponent of $S$. It has the following form (see \cite[III.1]{Be})
\[
    \phi(\lambda)=b\lambda+\int_{(0,\infty)}(1-e^{-\lambda t})\mu(dt)\,.
\]
Here, $b\geq 0$ is called the drift of $S$ and $\mu$ is a measure on $(0,\infty)$ satisfying $\int_{(0,\infty)}(1\wedge t)\mu(dt)<\infty$ called the L\' evy measure of $S$\,.

Let $B=(B_t,\P_x)_{t\geq 0,\, x\in \R^d}$ be the Brownian motion in $\R^d$ $(d\geq 1)$ independent of the subordinator $S$. We define the subordinate Brownian motion $X=(X_t,\P_x)_{t\geq 0,\, x\in \R^d}$ by $X_t=B_{S_t}$, $t\geq 0$. It is a L\' evy process (see \cite[Theorem 30.1]{S}) such that  
\[
    \E e^{i\xi\cdot X_t}=e^{-t\phi(|\xi|^2)},\quad \xi\in \R^d\,.
\]  
Moreover, it has transition density $p(t,x,y)=p(t,y-x)$ and it is of the form
\[
    p(t,x)=\int_{(0,\infty)}(4\pi s)^{-d/2}e^{-\frac{|x|^2}{4s}}\P(S_t\in ds)\,.
\]

Taking $\phi(\lambda)=\lambda^{\alpha/2}$ with $\alpha\in (0,2)$, we obtain rotationally invariant symmetric $\alpha$-stable process $X$. Its infinitesimal generator is the fractional Laplacian $-(-\Delta)^{\alpha/2}$ defined by
\[
    -(-\Delta)^{\alpha/2}u(x)=\int_{\R^d\setminus\{0\}}(u(x+y)-u(x)-\nabla u(x)\cdot y 1_{B_1(0)}(y))\frac{c_{d,\alpha}}{|y|^{d+\alpha}}\,dy
\]
and the following heat kernel estimate holds (see \cite{BG})
\begin{equation}\label{eq:hke-stable}
    p^{(\alpha)}(t,x,y)\asymp t^{-d/\alpha}\wedge \frac{t}{|x-y|^{d+\alpha}}=\phi^{-1}(t^{-1})^{d/2}\wedge \frac{t\phi(|x-y|^{-2})}{|x-y|^d}\,.
\end{equation}

Note that in this case heat kernel estimate can be expressed just in terms of the Laplace exponent and its inverse function. Such type of estimate will continue to hold if the scaling of the process $X$ is strictly between $0$ and $2$ (see \cite{BGR} and Corollary \ref{cor:classical}), but if the scaling fails to satisfy this condition we will see that  different form of heat kernel estimates appear. The approach in this paper is more general and it will essentially hold for subordinate Brownian motions with scaling order that is strictly greater than $0$.

We introduce the following scaling conditions for a function $f\colon (0,\infty)\rightarrow (0,\infty)$\,\,:
\begin{itemize}

\item[{\bf (L)}] there exist $\gamma>0$, $\lambda_L\geq 0$ and $C_L>0$ such that 
\begin{equation}\label{eq:lowerscaling}
    \frac{f(\lambda x)}{f(\lambda)}\geq C_Lx^{\gamma}\qquad \text{ for all }\quad \lambda>\lambda_L\ \text{ and }\ x\geq 1\,,
\end{equation}
\item[{\bf (U)}] there exist $\delta>0$, $\lambda_U\geq 0$ and $C_{U}>0$ such that 
\begin{equation}\label{eq:upperscaling}
    \frac{f(\lambda x)}{f(\lambda)}\leq C_{U}x^{\delta}\qquad \text{ for all }\quad \lambda> \lambda_U\ \text{ and }\ x\geq 1\,.
\end{equation}
\end{itemize}

If $f$ is non-decreasing, then {\bf (L)} and {\bf (U)} are actually doubling conditions meaning that it is enough that they hold for some $x>1$. 

Scaling conditions will be interesting for the Laplace exponent $\phi$ of a subordinator $S$ and for the function $H:(0,\infty)\rightarrow [0,\infty)$ defined by $H(\lambda):=\phi(\lambda)-\lambda\phi'(\lambda)$. The function $H$ appeared in the work of Jain and Pruitt \cite{JP}, where, in particular,  asymptotic properties of lower (upper) tail probabilities of subordinators were studied. We will need estimates of the tail probabilities of subordinators and in this sense tail estimates in Section \ref{sec:subord} represent an upper tail counterpoint of the results from \cite{JP} (see Proposition \ref{prop:subord_upper} and Proposition \ref{prop:subord-lower})\,. However, our method involves different techniques\,. 
 
It is not hard to show that  if $H$ satisfies {\bf (L)} or {\bf (U)}, respectively, then the same holds for $\phi$ in the case of a zero-drift subordinator (see Lemma \ref{lem:bf}). 
If $\phi$ satisfies {\bf (L)},  then processes with scaling order $0$  cannot be considered, e.g. geometric stable processes. These are the processes obtained by subordinating Brownian motion by geometric stable subordinators, that is the subordinators with the Laplace exponent of the form $\phi(\lambda)=\log(1+\lambda^{\beta/2})$ with $\beta\in (0,2]$. The corresponding near diagonal estimate is infinite for $d>\beta$. Indeed, by \cite[(3.8)]{SSV2}, 
\[
    p(1,x)=\int_0^\infty e^{-t} p^{(\beta)}(t,x)\,dt\geq c_1 \int_0^{|x|^\beta} e^{-t} t|x|^{-d-\beta}\,dt\geq   c_2|x|^{-d+\beta},
\]
hence $\lim\limits_{|x|\to 0+}p(1,x)=+\infty\,.$ Thus  estimating heat kernel in this case makes sense only away from the diagonal.

Actually, the lower scaling property of $\phi$  is  equivalent to the near diagonal upper bound of the heat kernel and this is the first result of this paper. 

\begin{Thm}\label{tm:upper-equiv}
Let $S=(S_t)_{t\geq 0}$ be a subordinator 
  and let $X=(X_t)_{t\geq 0}$ be the corresponding subordinate Brownian motion in $\R^d$ with the transition density $p(t,x,y)=p(t,y-x)$. Then there exist $C>0$ and $\lambda_L\geq 0$ such that 
  \begin{equation}\label{eq:on-up}
      p(t,x)\leq C\phi^{-1}(t^{-1})^{d/2}\quad \text{ for all }\quad 0<t<\phi(\lambda_L)^{-1}\,\, \text{ and }\,\, x\in \R^d
  \end{equation}
  if and only if $\phi$ satisfies {\bf (L)}.
\end{Thm}



The main result of the paper is the following sharp heat kernel estimate. 

\begin{Thm}\label{tm:main}
 Let $S=(S_t)_{t\geq 0}$ be a subordinator with zero drift 
  and let $X=(X_t)_{t\geq 0}$ be the corresponding subordinate Brownian motion in $\R^d$ with the transition density $p(t,x,y)=p(t,y-x)$\,.
 \begin{itemize}
  \item[(i)] If $\phi$ satisfies {\bf (L)}, then there exists a constant $\kappa\in (0,1)$ such that for all $0<t<\kappa\phi(\lambda_L)^{-1}$ and $x\in \R^d$ satisfying $t\phi(|x|^{-2})\geq 1$, the following  near diagonal estimate holds
  \[
   p(t,x)\asymp \phi^{-1}(t^{-1})^{d/2}\,.
  \]
  \item[(ii)] If $H$ satisfies {\bf (L)} and {\bf (U)} with $\delta<2$,  then there exist  constants 
    $\kappa, \eta, \theta\in (0,1)$, $C\geq 1$ and  $a_L,a_U>0$ such that for all $0<t<\kappa\phi(\lambda_L)^{-1}$ and $x\in \R^d$ satisfying $|x|< \theta\lambda_L^{-1/2}\wedge \eta\lambda_U^{-1/2}$ and $t\phi(|x|^{-2})\leq 1$, 
  the following off-diagonal estimates hold
 \begin{align*}
 p(t,x)&\leq C \left(t|x|^{-d}H(|x|^{-2})\vee \phi^{-1}(t^{-1})^{d/2}e^{-a_U|x|^2\phi^{-1}(t^{-1})}\right) \\
  p(t,x)&\geq C^{-1}\left(t|x|^{-d}H(|x|^{-2})\vee \phi^{-1}(t^{-1})^{d/2}e^{-a_L|x|^2\phi^{-1}(t^{-1})}\right) \\
 \end{align*}
 \end{itemize}
\end{Thm}
Note that, in the case $\lambda_L=\lambda_U=0$,  the heat kernel estimates in Theorem \ref{tm:main} are global in space and time. Since $t\phi(|x|^{-2})\leq 1$ is equivalent to $|x|^2\phi^{-1}(t^{-1})\geq 1$, we get the following off diagonal estimates too:
 \begin{align*}
 p(t,x)&\leq C \left(t|x|^{-d}H(|x|^{-2})\vee|x|^{-d}e^{-\tilde{a}_U|x|^2\phi^{-1}(t^{-1})}\right) \\
  p(t,x)&\geq C^{-1}\left(t|x|^{-d}H(|x|^{-2})\vee|x|^{-d}e^{-\tilde{a}_L|x|^2\phi^{-1}(t^{-1})}\right)\,, \\
 \end{align*}
for some $\tilde{a}_U>0$ and $\tilde{a}_L=a_L$.

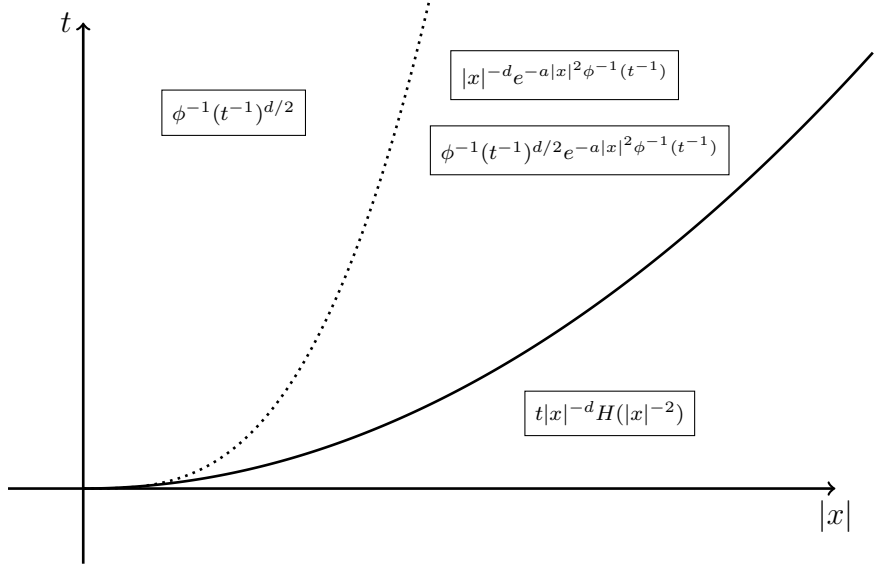
\begin{figure}\label{fig:1}
\begin{tikzpicture}
  \draw[->,line width=1pt] (-5,-2) -- (6,-2) node[below] {$|x|$};
  \draw[->,line width=1pt] (-4,-3) -- (-4,4.2) node[left] {$t$};
\draw[scale=1,samples=200,domain=-4:6.5,smooth,variable=\x, line width=1pt] plot ({\x},{(\x+4)^2/19-2});
\draw[scale=1,samples=200,domain=-3.998:0.6,smooth,variable=\x,dotted, line width=1pt] plot ({\x},{exp(3*ln(\x+4))/15-2});
\node[draw] at (-2,3) {\tiny $\phi^{-1}(t^{-1})^{d/2}$};
\node[draw] at (2.4,3.5) {\tiny $|x|^{-d}e^{-a|x|^2\phi^{-1}(t^{-1})}$};
\node[draw] at (2.6,2.5) {\tiny $ \phi^{-1}(t^{-1})^{d/2}e^{-a|x|^2\phi^{-1}(t^{-1})}$};
\node[draw] at (3,-1) {\tiny $t|x|^{-d}H(|x|^{-2})$};
\end{tikzpicture}
\caption{Regions of heat kernel estimates. Dotted line corresponds to  $t\phi(|x|^{-2})=1$, while  'full' line corresponds to $tH(|x|^{-2})=e^{-a|x|^2\phi^{-1}(t^{-1})}$.}
\end{figure}

A novelty in this result is the appearance of the term of the form \[\phi^{-1}(t^{-1})^{d/2}e^{-a|x|^2\phi^{-1}(t^{-1})}\,\,( \text{or }\,\,|x|^{-d}e^{-a|x|^2\phi^{-1}(t^{-1})})\] in the off-diagonal estimate. It could be explained as an intermediate term between on-diagonal estimate and the classical off-diagonal estimate $t|x|^{-d}H(|x|^{-2})$ that involves tail estimate of the L\' evy measure of the subordinate Brownian motion $X$ (see Figure \ref{fig:1})\,.

A natural question that arises is in what situations the functions $H$ and $\phi$ are comparable. 
It turns out that this is equivalent to the property that $\phi$ satisfies {\bf(U)} with $\delta<1$ (see Proposition \ref{prop:comp}).
In other words, $H$ can be replaced by $\phi$ in Theorem \ref{tm:main} precisely when the scaling is strictly between $0$ and $2$ and in this case the estimate from Theorem \ref{tm:main} reduces to the estimate that is possible to deduce from results already known from previous works (see e.g. \cite[eq. (2)]{BGR}). 
\begin{Cor}\label{cor:classical}
 Let $S=(S_t)_{t\geq 0}$ be a subordinator with the Laplace exponent $\phi$ satisfying {\bf (L)} and {\bf (U)} with $\delta<1$ and let $X=(X_t)_{t\geq 0}$ be the corresponding subordinate Brownian motion in $\R^d$ with the transition density $p(t,x,y)=p(t,y-x)$\,. Then 
 there exist constants $\kappa, \eta, \theta \in (0,1)$ such that for all $0<t<\kappa\phi(\lambda_L)^{-1}$ and $x\in \R^d$ satisfying $|x|< \theta\lambda_L^{-1/2}\wedge \eta\lambda_U^{-1/2}$, 
 \[
  p(t,x)\asymp \phi^{-1}(t^{-1})^{d/2}\wedge t|x|^{-d}\phi(|x|^{-2})\,.
 \]
\end{Cor}

In some cases it is easier to work with the second derivative of $\phi$ than $H$. 
\begin{Prop}\label{prop:second_der}
If  $\phi$ is the Laplace exponent of a zero drift subordinator such that $H$ satisfies {\bf (L)}, then there exists $M\geq 1$ such that 
\[
H(\lambda)\asymp \lambda^2(-\phi''(\lambda)),\quad \lambda>M\lambda_L\,.
\]
\end{Prop}
 A similar result has already appeared (see e.g. \cite[Lemma 5.1]{JP}). 

\begin{Rem}  Concerning comparability of certain function involving $\phi$ there are more interesting results. It turns out that if $H$ satisfies {\bf (L)} and {\bf (U)} with $\delta<2$, then the functions  $\phi(\lambda)$ and $\lambda\phi'(\lambda)$ are comparable. Conversely, comparability of these two functions implies that  $\phi$ satisfies {\bf(L)} (see Proposition \ref{prop:compL}).
\end{Rem}

Let us return to the new features of our result and consider examples with order of scaling that is not strictly between $0$ and $2$\,.
\begin{Ex} \label{ex:1}
Let $\phi(\lambda)=\frac{\lambda}{\log(1+\lambda^{\beta/2})}$, where $\beta\in (0,2)$\,. Here,
\[
 \phi^{-1}(\lambda)\asymp \begin{cases}
                           \lambda^{\frac{2}{2-\beta}} & 0<\lambda<2\\
                           \lambda\log\lambda & \lambda\geq 2
                          \end{cases}
                          \qquad \qquad
 H(\lambda)\asymp \begin{cases}
                   \lambda^{1-\beta/2} & 0<\lambda<2\\
                   \frac{\lambda}{(\log\lambda)^2} & \lambda\geq 2\,.
                  \end{cases}
\]
Hence, $H$ satisfies {\bf (L)} and {\bf (U)} with $\lambda_L=\lambda_U=0$ and $\delta<2$, but $\phi$ and $H$ are not comparable for large values of $\lambda$ (compare with Proposition \ref{prop:comp}).  For example, if $0<t<1/2$ and $0<|x|<1/2$ satisfy $t\phi(|x|^{-2})\leq 1$, then we get the off-diagonal estimate
\[
 p(t,x)\asymp \frac{t}{|x|^{d+2}(\log\frac{1}{|x|})^2}\vee t^{-d/2}\left(\log\frac{1}{t}\right)^{-d/2}e^{-a\frac{|x|^2}{t}\log\frac{1}{t}}\,,
\]
where $a=a_L$ or $a=a_U$ dependeing whether we consider lower or upper bound.  
For other $t$ and $x$ we can express estimates explicitly similarly by using behavior of $H$ and $\phi^{-1}$. In particular, for $t>1/2$ and $|x|>1/2$ we get stable-like estimates
\[
p(t,x)\asymp t^{-d/(2-\beta)}\wedge \frac{t}{|x|^{d+2-\beta}}\,.
\]
The first estimate gives an answer how a sharp upper bound in  \cite{Mi3} for a L\' evy processes should look like. In \cite[Theorem 4]{KSz} a similar upper estimate for this process has been obtained. 
\end{Ex}

\begin{Ex}\label{ex:2}
 Let $\phi(\lambda)=\frac{\lambda}{\log(1+\lambda)}-1$ be the Laplace exponent of the conjugate gamma subordinator without killing. In this case $\phi(\lambda)\sim \frac{\lambda}{2}$ and $H(\lambda)\sim\frac{\lambda^2}{12}$ as $\lambda \to 0+$. Hence, 
 \[
 \phi^{-1}(\lambda)\asymp \begin{cases}
                           \lambda& 0<\lambda<2\\
                           \lambda\log\lambda & \lambda\geq 2
                          \end{cases}
                          \qquad \qquad
 H(\lambda)\asymp \begin{cases}
                   \lambda^2& 0<\lambda<2\\
                   \frac{\lambda}{(\log\lambda)^2} & \lambda\geq 2\,.
                  \end{cases}
\]
Therefore, $\phi$ and $H$ are not comparable at all, but $H$ satisfies {\bf (L)} and {\bf (U)} with $\lambda_L=0$ and $\lambda_U=2$. In this case, we get the estimate only for small values of $t$ and $|x|$, which are basically the same as the corresponding estimates in Example \ref{ex:1}.
\end{Ex}


A crucial step in estimating heat kernel are upper tail estimates of subordinators (see Proposition \ref{prop:subord_upper}, Proposition \ref{prop:subord-low} and  Lemma \ref{lem:subord-tail}). Put together, if $H$ satisfies {\bf (U)} with $\delta<2$, we obtain the following estimate
\begin{equation}\label{eq:subord-tailestimate}
    \P(S_t\geq r)\asymp t H(r^{-1})
\end{equation}
for $t>0$ and $0<r<M\lambda_U^{-1}$ satisfying $t\phi(r^{-1})\leq e^{-1}(1-\varepsilon)$, where $M,\varepsilon\in (0,1)$ and constants in (\ref{eq:subord-tailestimate}) depend on $\varepsilon$\,.
As already mentioned, estimate (\ref{eq:subord-tailestimate}) could be understood as a uniform version of \cite[Theorem 5.1]{JP} for upper tails of subordinators.

As an application, estimates of heat kernel obtained in Theorem \ref{tm:main} can be used to estimate Green function (or potential) in the case of transient subordinate Brownian motion $X=(X_t)_{t\geq 0}$. To be more precise, we say that $X$ is transient if $\P_0(\lim\limits_{t\to\infty}|X_t|=\infty)=1$\,. A necessary and sufficient condition for transience is that there exists $a>0$ so that $\int_0^a \frac{y^{d-1}}{\phi(y^2)}\,dy<\infty$ (see \cite[Corollary 37.6]{S})\,. If $X$ is transient, then (see \cite[Theorem 35.4]{S}) the following function is well-defined and finite 
\[
 G(x,y)=G(x-y)=\int_0^\infty p(t,x-y)\,dt, \quad x,y\in \R^d,\ x\not=y
\]
and it is called the Green function (or potential) of $X$\,.

The following result gives behavior of the Green function when the scaling order is strictly greater than $0$. 
Such estimates are usually proved by using different techniques (e.g. using potential measure of the subordinator), but since we have obtained sharp heat kernel estimates it is almost a direct consequence.  

\begin{Cor}\label{cor:green}
 Let $S=(S_t)_{t\geq 0}$ be a 
 subordinator with zero drift such that its Laplace exponent $\phi$ satisfies {\bf (L)} with $\lambda_L=0$
  and let $X=(X_t)_{t\geq 0}$ be the corresponding subordinate Brownian motion in $\R^d$. Assume that $X$ is transient and that in the case $d\leq 2$ the Laplace exponent $\phi$ satisfies additionally {\bf (U)} with $\delta<d/2$ and $\lambda_U=0$. 
  Then
  \[
   G(x)\asymp \frac{1}{|x|^d\phi(|x|^{-2})},\quad x\in \R^d\,.
  \]
\end{Cor}
 

Green function of the conjugate gamma process from Example \ref{ex:2}, for $d\geq 3$, can be estimated as
\begin{equation}\label{eq:green-conjgamma}
    G(x)\asymp \begin{cases}
    |x|^{2-d}\,\log\tfrac{1}{|x|} & |x|< \tfrac{1}{2}\\
     |x|^{2-d} & |x|\geq \tfrac{1}{2}\,.\\
 \end{cases}
\end{equation}
This shows how close this process is to the Brownian motion, since the Green function of the Brownian motion in $\R^d$ ($d\geq 3$) is the Newtonian potential $G^{(2)}(x)=c|x|^{2-d}$\,.

Estimate in Corollary \ref{cor:green} for subordinate Brownian motions is known in many cases when the scaling order is strictly between $0$ and $2$ (see \cite{SV2, KSV3}). For some cases it has been calculated when the scaling order is allowed to be near $2$ (see \cite{Mi2}) and recently it has been proved that the condition {\bf (L)} for $\phi$ is actually equivalent to the lower bound of the Green function in Corollary \ref{cor:green} for $|x|$ small (see \cite{G}).
The form of the Green function behavior in the case of zero order scaling is different (for details, see \cite{SSV2, KM,Mi3})\,.

The structure of the paper is as follows. In Section \ref{sec:subord} tail estimates of subordinators are obtained. Further, some properties of functions $\phi$ and  $H$ and their relationships are proved. Some of the results in this section could be also considered to be of independent interest. Section \ref{sec:hke} is devoted to heat kernel estimates. It consists of three subsections, first two for upper and lower estimates and in the last one we prove Green function estimates and the main results.

\section{Subordinators and their tail estimates}
\label{sec:subord}

Let $S=(S_t)_{t\geq 0}$ be a subordinator. Recall that the Laplace transform of $S_t$ is given by
\[
	\E [e^{-\lambda S_t}]=e^{-t\phi(\lambda)}\qquad \phi(\lambda)=b\lambda+\int_{(0,\infty)} (1-e^{-\lambda y})\mu(dy),\ \lambda>0,
\]
where $\phi$ is called the Laplace exponent, $b\geq 0$ is the drift and $\mu$ is the  L\' evy measure of $S$ meaning that $\mu$ is a measure on $(0,\infty)$ satisfying $\int_{(0,\infty)}(1\wedge y)\mu(dy)<\infty$\,. 

The Laplace exponent $\phi$ belongs to the class of Bernstein functions $\mathcal{BF}=\{f\in C^\infty(0,\infty)\colon f\geq 0,\,  (-1)^{n-1}f^{(n)}\geq 0,\ n\in \N\}\,$.
It is known that every Bernstein function $f$ has a unique representation 
\begin{equation}\label{eq:bf-repr}
	f(\lambda)=a+b\lambda+\int_{(0,\infty)} (1-e^{-\lambda y})\mu(dy),
\end{equation}
where $a,b\geq 0$ and $\mu$ is a L\' evy measure (cf. \cite[Chapter 3]{SSV})\,.

The following formula that follows from Fubini theorem will be useful
\begin{equation}\label{eq:bf-tail}
\frac{f(\lambda)}{\lambda}=\frac{a}{\lambda}+b+\int_{(0,\infty)}e^{-\lambda t}\mu(t,\infty)\,dt,\quad \lambda>0\,.
\end{equation}
Further, the following inequality will be often used:
\begin{align}\nonumber
    \lambda f'(\lambda)&=\lambda b+\int_{(0,\infty)} \lambda  te^{-\lambda t}\mu(dt)\\&\leq\lambda b+ \int_{(0,\infty)}(1-e^{-\lambda t})\mu(dt)\leq f(\lambda),\quad \lambda>0\,. \label{eq:tmp3412}
\end{align}
A direct consequence of inequality (\ref{eq:tmp3412}) is non-negativity of the function $H$\,.

The first lemma gives some general properties of Bernstein function needed in the sequel. Although most of results are known, we present some alternative proofs. 

\begin{Lem}\label{lem:bf}
Let $\phi\in \mathcal{BF}$. 
\begin{itemize}
\item[(a)] For any $\lambda>0$ and $x\geq 1$,
\[
    \phi(\lambda x)\leq x\phi(\lambda)\qquad \text{ and }\qquad H(\lambda x)\leq x^2H(\lambda)\,.
\]
\item[(b)] 
Assume that the drift of $\phi$ in the representation (\ref{eq:bf-repr}) is $b=0$. If $H$ satisfies {\bf (L)} (resp. {\bf (U)}), then $\phi$ satisfies {\bf (L)}(resp. {\bf (U)}).  
\end{itemize}
\end{Lem}
\proof
(a) Without loss of generality we may assume that $a=b=0$ in the represtntation (\ref{eq:bf-repr}) of $\phi$. The result concerning $\phi$ is already known and follows from concavity of $\phi$. Here we present new proofs of both inequalities that follow from the following elementary inequalities
\begin{align*}
    &1-e^{-\lambda x}\leq x(1-e^{-\lambda})\,,\\ &1-e^{-\lambda x}-\lambda x e^{-\lambda x}\leq x^2(1-e^{-\lambda}-\lambda e^{-\lambda}),\quad  x\geq 1,\ \lambda>0\,.
\end{align*}
Indeed, 
\begin{align*}
\phi(\lambda x)&=\int_{(0,\infty)}(1-e^{-\lambda x t})\mu(dt)\leq x\int_{(0,\infty)}(1-e^{-\lambda  t})\mu(dt)=x\phi(\lambda)
\end{align*}
and
\begin{align*}
H(\lambda x)&=\int_{(0,\infty)}(1-e^{-\lambda x t}-\lambda x te^{-\lambda x t})\mu(dt)\\&\leq x^2\int_{(0,\infty)}(1-e^{-\lambda  t}-\lambda t e^{-\lambda t})\mu(dt)=x^2H(\lambda)
\end{align*}
(b) Assume now that $H$ satisfies {\bf(U)} and that the drift is $b=0$. Then $\lim\limits_{\lambda\to\infty }\frac{\phi(\lambda)}{\lambda}=b=0$ (see \cite[Proposition I.2 (ii)]{Be}); hence
\[
\frac{\phi(\lambda)}{\lambda}=\int_\lambda^\infty \frac{H(s)}{s^2}\,ds\,.
\]
 By change of variable,
 \begin{align*}
 \phi(\lambda x)&=\lambda x\int_{\lambda x}^\infty \frac{H(s)}{s^2}\,ds=\lambda\int_\lambda^\infty \frac{H(sx)}{s^2}\,ds\\
 &\leq C_U \lambda x^\delta \int_\lambda^\infty \frac{H(s)}{s^2}\,ds=C_Ux^\delta \phi(\lambda)\,.
 \end{align*}
Similar argument applies to {\bf(L)}.
\qed

Note that Lemma \ref{lem:bf} (a) suggests that $\phi$ and $H$ could have different scaling properties and this can be seen in Example \ref{ex:2}. To be more precise, for $\phi(\lambda)=\frac{\lambda}{\log(1+\lambda)}-1$ the bounds from Lemma \ref{lem:bf} (a) are attained:
\[
    \lim\limits_{\lambda\to 0+}\frac{\phi(\lambda x)}{\phi(\lambda)}=x\qquad \text{ and }\qquad \lim\limits_{\lambda\to 0+}\frac{H(\lambda x)}{H(\lambda)}=x^2\,.
\]

The following tail estimate of the subordinator reveals a probabilistic connection between the functions $H$ and $\phi$ and it will play an important role in obtaining upper off-diagonal estimates of the heat kernel. Before this we need an auxiliary lemma.

\begin{Lem}\label{lem:levy}
	Let $\mu$ be the L\' evy measure of a subordinator with the Laplace exponent $\phi$ and the L\' evy measure $\mu$. Then, for any $r>0$, 
	\begin{align*}
		&\text{(i)}\quad \mu(r,\infty)+r^{-2}\int_{(0,r]}y^2\mu(dy)=2r^{-2}\int_0^ry\mu(y,\infty)\,dy\\
		&\text{(ii)}\quad \int_0^ry\mu(y,\infty)\,dy\leq er^2 H(r^{-1})\\
		&\text{(iii)}\quad \int_{(0,r]}y\mu(dy)\leq er\phi(r^{-1})\,.
	\end{align*}
\end{Lem}
\proof
	(i) By Fubini theorem we get
	\begin{align*}
	    \int_{(0,r]}y^2\mu(dy)&=\int_{(0,r]}\int_0^y 2z\,dz\mu(dy)=\int_0^r 2z\mu((z,r])\,dz\\
	    &=\int_0^r 2z(\mu(z,\infty)-\mu(r,\infty))\,dz=2\int_0^r z\mu(z,\infty)\,dz-r^2\mu(r,\infty)\,.
	\end{align*}
	(ii) Taking derivative in 
	\[
		\frac{\phi(\lambda)}{\lambda}=b+\int_0^\infty e^{-\lambda y}\mu(y,\infty)\,dy
	\]
	it follows that 
	\begin{align}\label{eq:lem-subord}
		\frac{H(\lambda)}{\lambda^2}=\frac{\phi(\lambda)-\lambda\phi'(\lambda)}{\lambda^2}=\int_0^\infty e^{-\lambda y}y\mu(y,\infty)\,dy\geq e^{-1}\int_0^{\lambda^{-1}}y\mu(y,\infty)\,dy\,.
	\end{align}
	The desired inequality is obtained by choosing $\lambda=r^{-1}$\,.\\
	(iii) Since $1-e^{-x}\geq xe^{-x}$ for $x\geq 0$, 
	\[
		\phi(\lambda)\geq \int_{(0,\infty)}\lambda ye^{-\lambda y}\mu(dy)\geq e^{-1} \lambda \int_{(0,\lambda^{-1}]}y\mu(dy)
	\]
	and now it is enough to take $\lambda=r^{-1}$\,.\\
\qed	

\begin{Prop}\label{prop:subord_upper}
	Let $S=(S_t)_{t\geq 0}$ be a subordinator with the Laplace exponent $\phi$. Then there exists a constant $C_S>0$ such that 
		\[
		\P(S_t\geq r(1+et\phi(r^{-1})))\leq  
		C_StH(r^{-1})
	\]
	for all $t,r>0$\,. In particular, if $0<\varepsilon<1$, then
	\[
	    \P(S_t\geq r)\leq C_S t H(\varepsilon^{-1}r^{-1}) \leq C_S\varepsilon^{-2} tH(r^{-1})
	\]
	for all $t,r>0$ satisfying  $t\phi(r^{-1})\leq e^{-1}(1-\varepsilon)$.\end{Prop}
\proof
	Let 
	\[
		\wt{S}_t:=S_t-t\int_{(0,r]}y\mu(dy),\quad t\geq 0\, , r>0\,.
	\]
	Then 
	\begin{align*}
		\E[e^{i\theta \wt{S}_t}]&=\exp{\left\{-t\phi(-i\theta)-i\theta t\int_{(0,r]}y\mu(dy)\right\}}\\&=\exp{\left\{-t\int_{(0,\infty)}\left(1-e^{i\theta y}+i\theta y1_{\{|y|\leq r\}}\right)\mu(dy)\right\}}
	\end{align*}
	which shows that $\{\wt{S}_t\}_{t\geq 0}$ is a L\' evy process in $\R$ with the L\' evy exponent $\psi(\theta)=\int_{(0,\infty)}(1-e^{i\theta y}+i\theta y1_{\{|y|\leq r\}})\mu(dy)$\,.
	Hence it follows that for any $f\in C_0^2(\R):=\{f\in C^2(\R)\colon f,f',f''\in C_0(\R)\}$, where $C_0(\R)=\{f\in C(\R):\lim\limits_{x\to\pm \infty}f(x)=0\}$, the infinitesimal generator of $\wt{S}$ is given by (see \cite[Theorem 31.5]{S})
	\[
		\mathcal{L}f(x)=\int_{(0,\infty)}\left(f(x+y)-f(x)-f'(x)y1_{\{|y|\leq r\}}\right)\mu(dy)\,.
	\]
	Let $g\in C^2(\R)$ with bounded first and second derivatives such that 
	\begin{align*}
	  &0\le g(x)\le 1 \text{ for all }x\in \R,\quad    g(x)=0\text{ for } x<\frac{1}{4}\quad \text{ and }\quad  g(x)=1\text{ for } x\geq \frac{3}{4}\,
	\end{align*}
	and define, for $n\ge 2$,
	\[
	    f_n(x)=\begin{cases}
                g(x) & x\leq 1\\
                1 & 1< x< n\\
                g(n+1-x) & x\geq n\,.
                \end{cases}
	\]
	Then $f_n\in C_0^2(\R)$ and 
	\[
	    \|f_n\|_\infty=1\,,\|f_n'\|_\infty\leq \|g'\|_\infty \quad \text{ and }\quad \|f_n''\|_\infty\leq \|g''\|_\infty \quad \text{ for all }\,\,n\in \N. 
	\]
Let $n\in \N$ and set $f_{n,r}(y):=f_n(\frac{y}{r})$\,. Then
	\begin{align*}
		f_{n,r}(x+y)-f_{n,r}(x)-\frac{d}{dx}f_{n,r}(x)y1_{\{|y|\leq r\}}&=f_n\left(\tfrac{x+y}{r}\right)-f_n\left(\tfrac{x}{r}\right)-f_n'\left(\tfrac{x}{r}\right)\tfrac{y}{r}1_{\{|\frac{y}{r}|\leq 1\}}\\&\leq \tfrac{1}{2}\|f_n''\|_\infty \left(\tfrac{y}{r}\right)^21_{\{|y|\leq r\}}+\|f_n\|_\infty1_{\{|y|>r\}}
	\end{align*}
	and Lemma \ref{lem:levy} yields
	\begin{align*}
		\mathcal Lf_{n,r}(x)&=\int_{(0,\infty)}\left(f_n\left(\tfrac{x+y}{r}\right)-f_n\left(\tfrac{x}{r}\right)-f'\left(\tfrac{x}{r}\right)\tfrac{y}{r}1_{\{|\frac{y}{r}|\leq 1\}}\right)\mu(dy)\\
		&\leq  \tfrac{1}{2}\|f_n''\|_\infty r^{-2}\int_{(0,r]}y^2\mu(dy)+\|f_n\|_\infty \int_{(r,\infty)}\mu(dy)\\
		&\leq 2(\|g''\|_\infty\vee 1)r^{-2}\int_0^r y\mu(y,\infty)\,dy\\
		&= C_SH(r^{-1})
	\end{align*}
	with $C_S:= 2e(\|g''\|_\infty\vee 1)$\,.
	By Dynkin formula (see \cite[Proposition VII.1.2 and Proposition VII.1.6]{RY}) and the last display it follows that 
	\begin{align*}
		  \E[ f_{n,r}(\wt{S}_t)]&=\E\int_0^t \mathcal{L}f_{n,r}(\wt{S}_s)\,ds
		\leq C_StH(r^{-1})\,.
	\end{align*}
	Since $\uparrow\lim\limits_{n\to\infty}f_{n,r}(x)=g_r(x):=g(\frac{x}{r})$  
	we may use monotone convergence theorem to conclude that 
\[
     \E[ g_r(\wt{S}_t)]=\lim_{n\to\infty}\E[ f_{n,r}(\wt{S}_t)]
		\leq C_StH(r^{-1})\,.
\]
	Hence, Lemma \ref{lem:levy} (iii) and the last display finally yield
	\begin{align*}
		\P(S_t\geq r(1+et\phi(r^{-1}))&\leq \P(\wt{S}_t\geq r)\leq \E[ g_r(\wt{S}_t)]\\&=\E\int_0^t \mathcal{L}g_r(\wt{S}_s)\,ds
		\leq C_StH(r^{-1})\,.
	\end{align*}
	Let $\theta >0$ and $t,r>0$ such that $t\phi(r^{-1})\leq \theta$\,, then    
\begin{equation}\label{eq:tmp_10023}
\P(S_t\geq r(1+e\theta))\leq  \P(S_t \geq r(1+et\phi(r^{-1}))\leq C_StH(r^{-1})\,.
\end{equation}
Now, if $t\phi(r^{-1})\leq \frac{\theta}{1+e\theta}$, by Lemma \ref{lem:bf} (a),
\[
    t\phi(r^{-1}(1+e\theta))\leq t\phi(r^{-1})(1+e\theta)\leq \theta,
\]
hence, by (\ref{eq:tmp_10023}) and Lemma \ref{lem:bf} (a),
\begin{align*}
    \P(S_t\geq r)&=\P(S_t\geq r(1+e\theta)^{-1}(1+e\theta))\leq C_St H(r^{-1}(1+e\theta))\\&\leq (1+e\theta)^2C_S tH(r^{-1})\,. 
\end{align*}
Now take $\varepsilon=\frac{1}{1+e\theta}$\,.
\qed

%
%
%
%

The following estimate will be needed to get the lower bound estimates of heat kernel.

\begin{Prop}\label{prop:sub-low}
Let $S=(S_t)_{t\ge 0}$ be a subordinator with Laplace exponent $\phi$.
For any $\alpha, \beta>0$ such that $\alpha \beta>1$ there exist constants $\rho\in (0,1)$ and $\tau>0$ such that 
\[
    \P\left(\frac{1}{\alpha\phi^{-1}(\beta^{-1})}\leq S_t\leq \frac{1}{\phi^{-1}(\rho t^{-1})}\right)\geq \tau\quad \text{ for all }\quad t>0\,.
\]
\end{Prop}
\proof
    Let $\alpha \beta>1$. 
    Since $y\mapsto 1-e^{-y}$ is strictly increasing and $y\mapsto e^{-y}$ is strictly decreasing, Markov inequality implies
    \begin{align*}
    \allowdisplaybreaks
    \P\left(\tfrac{1}{\alpha\phi^{-1}(\beta^{-1})}\right.&\left.\leq S_t\leq \tfrac{1}{\phi^{-1}(\rho t^{-1})}\right)
    \ge
    1-\P\left(S_t>\tfrac{1}{\phi^{-1}(\rho t^{-1})}\right)-\P\left(S_t<\tfrac{1}{\alpha\phi^{-1}(\beta t^{-1})}\right)\\
    &=1-\P\left(1-e^{-\phi^{-1}(\rho t^{-1})S_t}>1-e^{-1}\right)-\P\left(e^{-\phi^{-1}(\beta t^{-1})S_t}>e^{-1/\alpha}\right)\\
    &\geq 1-\frac{\E\left[1-e^{-\phi^{-1}(\rho t^{-1})S_t}\right]}{1-e^{-1}}-\frac{\E\left[e^{-\phi^{-1}(\beta t^{-1})S_t}\right]}{e^{-1/\alpha}}\\
    &=1-\frac{1-e^{-t\phi(\phi^{-1}(\rho t^{-1}))}}{1-e^{-1}}-\frac{e^{-t\phi(\phi^{-1}(\beta t^{-1}))}}{e^{-1/\alpha}}\\
    &=1-\frac{1-e^{-\rho}}{1-e^{-1}}-e^{1/\alpha-\beta}\,.
    \end{align*}
    Since $1-e^{1/\alpha-\beta}>0$ we can choose $\rho\in (0,1)$ small enough so that $\tau:=1-\frac{1-e^{-\rho}}{1-e^{-1}}-e^{1/\alpha-\beta}>0$\,.
\qed

To obtain the lower bound the following result will play an important role. 
\begin{Prop}\label{prop:subord-low}
  If $S$ is a subordinator with the L\' evy measure $\mu$ and the drift term $b\ge 0$, then for all $t,r>0$ we have
 \[
  \P(S_t\geq r)\geq 1-e^{-t\mu(r,\infty)}\,.
 \]
\end{Prop}
\proof Let $r>0$ and 
let $T=(T_t)_{t\geq 0}$ and $V=(V_t)_{t\geq 0}$ be two independent L\' evy processes with L\' evy measures
\[
\mu\cdot 1_{(0,r]}\quad \text{ and }\quad \mu\cdot 1_{(r,\infty)}\quad \text{respectively}\,,
\]
and $T$ having drift $b$. Then $S_t=T_t+V_t$ for $t\geq 0$ and $T$ and $V$ are both subordinators. Note that $V$ is a compound Poisson process (since its L\' evy measure is finite)\,. Then we have
\[
 \P(S_t\geq r)\geq \P(V_t\geq r)=1-e^{-t\mu(r,\infty)},
\]
since the event $\{V_t\geq r\}$ occurs precisely when the  the compound Poisson process $V$ jumps for the first time before time $t$\,.

\qed

Since the tail of the L\' evy measure appears in the estimate in the previous result, it would be useful to know its behavior.

\begin{Lem}\label{lem:subord-tail}
Let $\mu$ be the L\' evy measure of a subordinator with the Laplace exponent $\phi$. Then
\[
 \mu(r,\infty)\leq 2eH(r^{-1})\,, \quad r>0\,.
\]

If $H$ satisfies {\bf(U)} with $\delta<2$, then there exists $M\in (0,1)$ such that 
\[
	\mu(r,\infty)\geq \tfrac{M^2}{2} H(r^{-1}),\qquad 0<r<M\lambda_U^{-1}.
\]	
\end{Lem}
\begin{Rem}
 Note that for $\lambda_U=0$ the estimates hold for all $r>0$. If H satisfies {\bf(U)}  with $\delta<2$, we may write
 \[
     \mu(r,\infty)\asymp H(r^{-1}),\qquad 0<r<M\lambda_U^{-1}.
 \]
\end{Rem}

\proof
	Taking derivative in (see (\ref{eq:bf-tail}))
	\begin{equation}\label{eq:tail-phi}
	 \frac{\phi(\lambda)}{\lambda}=b+\int_0^\infty e^{-\lambda y}\mu(y,\infty)\,dy
	\end{equation}
	we get
	\[
	 \frac{H(\lambda)}{\lambda^2}=\int_0^\infty e^{-\lambda y}y\mu(y,\infty)\,dy\,, \quad \lambda>0
	\]
	and changing variable it follows that
	\begin{equation}\label{eq:laplace-low-tail}
	 H(r^{-1})=\int_0^\infty e^{-y} y\mu(r y,\infty)\,dy\,,\quad r>0\,.
	\end{equation}
	The upper estimate follows immediately, since it follows from (\ref{eq:laplace-low-tail}) that 
	\[
	 H(r^{-1})\geq \int_0^1 e^{-y} y \mu(ry,\infty)\,dy\geq \tfrac{e^{-1}}{2}\mu(r,\infty)\,,
	\]
	hence
	\begin{equation}\label{eq:levy-up}
	 \mu(r,\infty)\leq 2eH(r^{-1})\,, \quad r>0\,.
	\end{equation}

	To get the lower bound we start with (\ref{eq:laplace-low-tail}) and use the upper estimate (\ref{eq:levy-up}) and {\bf (U)} to get, for any $M\in (0,1)$ and $0<r<\lambda_U^{-1}$,  the following
	\begin{align*}
	 H(r^{-1})&=\int_0^Me^{-y}y\mu(ry,\infty)\,dy+\int_M^\infty e^{-y}y\mu(ry,\infty)\,dy\\
	 &\leq 2e\int_0^Me^{-y}yH(r^{-1}y^{-1})\,dy+\mu(Mr,\infty)\int_M^\infty e^{-y}y\,dy\\
	 &\leq 2eC_U(2-\delta)^{-1}M^{2-\delta}H(r^{-1})+(1+M)e^{-M}\mu(Mr,\infty)\,.
	\end{align*}
	Choosing $M\in (0,1)$ so that $2eC_U(2-\delta)^{-1}M^{2-\delta}\leq \frac{1}{2}$ and using Lemma \ref{lem:bf} (a), for $0<r<\lambda_U^{-1}$, we get
	\[
	 \mu(Mr,\infty)\geq \tfrac{e^M}{2(1+M)}H(r^{-1})\geq \tfrac{e^MM^{2}}{2(1+M)} H((Mr)^{-1})\,.
	\]
	It is enough to note that $\tfrac{e^MM^{2}}{2(1+M)}\geq \frac{M^2}{2}\,. $
\qed

\begin{Prop}\label{prop:subord-lower}
Let $S=(S_t)_{t\geq 0}$ be a zero-drift subordinator such that $H$ satisfies {\bf (L)} and {\bf (U)} with $\delta<2$. There exist constants $M\in (0,1)$, $L>1$ and $c_S>0$ such that for all $t>0$ and 
$0<r<L^{-1}\lambda_L^{-1}\wedge M\lambda_U^{-1}$ satisfying $t\phi(r^{-1})\leq 1$, 
\[
 \P(r\leq S_t\leq Lr)\geq c_St H(r^{-1})\,.
\]
\end{Prop}

\proof
  By Lemma \ref{lem:bf}, Proposition \ref{prop:subord_upper}, Proposition \ref{prop:subord-low} and Lemma \ref{lem:subord-tail}   it follows that there exists $M\in (0,1)$ such that, for 
   $r\in (0,L^{-1}\lambda_L^{-1}\wedge M\lambda_U^{-1})$,
  \begin{align*}
   \P(r\leq S_t\leq Lr)&\geq \P(S_t\geq r)-\P(S_t\geq Lr)\\
   &\geq 1-e^{-t\mu(r,\infty)}-4C_StH(L^{-1}r^{-1})\\
   &\geq t\mu(r,\infty)e^{-t\mu(r,\infty)}-4C_SC_L^{-1}L^{-\gamma}tH(r^{-1})\\
   &\geq \frac{M^2}{4}tH(r^{-1})e^{-2etH(r^{-1})}-4C_SC_L^{-1}L^{-\gamma}tH(r^{-1})\\
   &\geq \frac{M^2}{4}tH(r^{-1})e^{-2et\phi(r^{-1})}-4C_SC_L^{-1}L^{-\gamma}tH(r^{-1})\\
   &\geq tH(r^{-1})\left(\frac{M^2}{4}e^{-2e}-4C_SC_L^{-1}L^{-\gamma}\right)\\
  \end{align*}
In  the first inequality  
  we have used that  \[t\phi(L^{-1}r^{-1})\leq C_L^{-1}L^{-\gamma}t\phi(r^{-1})\leq C_L^{-1}L^{-\gamma}<\frac{1}{2e}\] for $L>1$ chosen large enough, since {\bf (L)} also holds for $\phi$  by Lemma \ref{lem:bf}.
  Again, choosing $L>1$ large enough so that 
  $\frac{M^2}{4}e^{-2e}-4C_SC_L^{-1}L^{-\gamma}>0$
  we obtain the desired bound.
\qed

Now it is possible to prove the result that establishes comparability of $H$ and $\lambda^2(-\phi''(\lambda))$. 

\proof[Proof of Proposition \ref{prop:second_der}]
 Here we use elementary inequalities
 \[
     \tfrac{1}{2}x^2e^{-x}\leq 1-e^{-x}-xe^{-x}\leq \tfrac{1}{2}x^2,\quad x>0\,.
 \]
 Then
\begin{align*}
    \tfrac{1}{2}\lambda^2(-\phi''(\lambda))&=\int_{(0,\infty)}\tfrac{1}{2}\lambda^2t^2e^{-\lambda t}\mu(dt)\\&\leq \int_{(0,\infty)}(1-e^{-\lambda t}-\lambda t e^{-\lambda t})\mu(dt)=H(\lambda),\quad \lambda>0\,.
\end{align*}
 For $M\in (0,1)$ we obtain
 \begin{align*}
 H(\lambda)&\leq \int_{(0,M^{-1}\lambda^{-1}]}(1-e^{-\lambda t}-\lambda t e^{-\lambda t})\mu(dt)+\mu(M^{-1}\lambda^{-1},\infty)\\
 &\leq \tfrac{1}{2} e^{M^{-1}}\int_{(0,M^{-1}\lambda^{-1}]}\lambda^2 t^2e^{-\lambda t}\mu(dt)+2eH(M\lambda)\\
 &\leq \tfrac{1}{2} e^{M^{-1}}\lambda^2 (-\phi''(\lambda))+2eC_L^{-1}M^\gamma H(\lambda),
 \end{align*}
where in the second inequality we have used Lemma \ref{lem:subord-tail} and in the last {\bf (L)}. Choosing $M\in (0,1)$ so that $2eC_L^{-1}M^\gamma=\frac{1}{2}$
we get 
\[
    H(\lambda)\leq e^{M^{-1}}\lambda^2(-\phi''(\lambda)),\quad \lambda>M^{-1}\lambda_L\,.
\]  
\qed

This section ends with results that explore comparability of the functions $H$, $\phi$ and $\lambda\phi'(\lambda)$\,.

\begin{Prop}\label{prop:comp}
The upper scaling condition {\bf (U)} for $\phi$ with $\delta<1$ holds if and only if
\[
 H(\lambda)\asymp \phi(\lambda)\,,\quad\lambda>\lambda_U\,.
\]
\end{Prop}
\proof Assume first that {\bf (U)} holds for $\phi$ with $\delta<1$. 
  Changing variables in (\ref{eq:tail-phi}) we obtain
  \[
   \phi(\lambda)=\int_0^\infty e^{-y}\mu(\lambda^{-1}y,\infty)\,dy,\quad \lambda>0\,.
  \]
  Since $H(\lambda)\leq \phi(\lambda)$, by Lemma \ref{lem:subord-tail},  for any $M\in (0,1)$ we have
  \begin{align*}
   \phi(\lambda)&=\int_0^M e^{-y}\mu(\lambda^{-1}y,\infty)\,dy+\int_M^\infty e^{-y}\mu(\lambda^{-1}y,\infty)\,dy\\
		&\leq 2e\int_0^M\phi(\lambda y^{-1})\,dy+\mu(M\lambda^{-1},\infty)\int_M^\infty e^{-y}\,dy \\
		&\leq 2e C_U\phi(\lambda)\int_0^M y^{-\delta}\,dy+2eH(M^{-1}\lambda)e^{-M}\\
		&\leq 2eC_U(1-\delta)^{-1}M^{1-\delta}\phi(\lambda)+2eM^{-2}H(\lambda),
  \end{align*}
  where in the last inequality we have used Lemma \ref{lem:bf} (a)\,. Since $1-\delta>0$, we can choose $M>0$ small enough so that $2eC_U(1-\delta)^{-1}M^{1-\delta}=\frac{1}{2}$ to get that \[\tfrac{M^2}{4e}\phi(\lambda)\leq H(\lambda)\leq \phi(\lambda)\ \text{ for all }\ \lambda>\lambda_U\,.\]
  
  Assume now that there exists $c\in (0,1)$ such that $\phi(\lambda)\geq H(\lambda)\geq c\phi(\lambda)$ for all $\lambda>\lambda_U$. This implies
  \[
   \frac{\phi'(t)}{\phi(t)}\leq \frac{1-c}{t}\,, \quad t>\lambda_U\,.
  \]
  Integrating this inequality over $[\lambda,\lambda x]$ for $\lambda>\lambda_U$ and $x\geq 1$ we obtain
  \[
   \log\frac{\phi(\lambda x)}{\phi(\lambda)}\leq (1-c)\log x
  \]
  which is {\bf (U)} with $\delta=1-c<1$\,.
\qed

The last result of the section explores comparability of functions and lower scaling conditions. 

\begin{Prop}\label{prop:compL}
The lower and upper scaling conditions {\bf (L)}  and {\bf (U)} for $H$ with $\delta<2$ imply
\begin{equation}\label{eq:comp3}
 \lambda\phi'(\lambda)\asymp \phi(\lambda)\,,\quad\lambda>\lambda_L\vee M^{-1}\lambda_U\,,
\end{equation}
for some constant $M\in (0,1)$\,.
Conversely, 
\[
    \lambda\phi'(\lambda)\asymp \phi(\lambda)\,,\quad\lambda>\lambda_L
\]
 implies the lower scaling condition {\bf (L)} for $\phi$\,.
\end{Prop}
\proof
  The upper bound follows from (\ref{eq:tmp3412}). 
 Assume that {\bf (L)} and {\bf (U)} hold for $H$ with $\delta<2$.  Then, for $L\geq 1$,  Lemma \ref{lem:subord-tail} yields
 \begin{align*}
  \lambda\phi'(\lambda)&\geq \int_{L^{-1}\lambda^{-1}}^{\lambda^{-1}}\lambda te^{-\lambda t}\mu(dt)\geq L^{-1}e^{-1}\left(\mu(L^{-1}\lambda^{-1},\infty)-\mu(\lambda^{-1},\infty)\right)\\
  &\geq L^{-1}e^{-1}\left(c_1H(\lambda L)-c_2H(\lambda)\right)
  =c_2 L^{-1}e^{-1}H(\lambda)\left(c_1c_2^{-1}\frac{H(\lambda L)}{H(\lambda)}-1\right)\\
    &\geq c_2 L^{-1} e^{-1}H(\lambda)\left(c_1c_2^{-1}C_L L^\gamma-1\right),\quad \lambda>M^{-1}\lambda_U
 \end{align*}
 for some $M\in (0,1)$\,.
Choosing $L\geq 1$ large enough so that $c_1c_2^{-1}C_L L^\gamma-1>0$ we get $\lambda\phi'(\lambda)\geq c_3H(\lambda)$ for all $\lambda>\lambda_L\vee M^{-1}\lambda_U$ implying
\[
 \lambda\phi'(\lambda)\geq \frac{c_3}{1+c_3}\phi(\lambda)\quad \text{ for all }\quad \lambda>\lambda_L\vee M^{-1}\lambda_U\,.
\]

Assume that there exists $c\in (0,1)$ so that  $\lambda\phi'(\lambda)\geq c\phi(\lambda)$ for all $\lambda>\lambda_L$\,. Similarly as in the proof of Proposition \ref{prop:comp} we obtain that for all $x\geq 1$ and $\lambda>\lambda_L$ the following inequality holds
\[
   \log\frac{\phi(\lambda x)}{\phi(\lambda)}\geq c\log x
  \]
and this is {\bf (L)} with $\gamma=c$\,.
\qed

\section{Heat kernel estimates for SBM}
\label{sec:hke}

%



In this section we obtain estimates of transition density of the subordinate Brownian motion. 

Let $S=(S_t)_{t\geq 0}$ be a subordinator defined on a probability space $(\Omega,\sig,\P)$ with the Laplace exponent $\phi$ and let 
$B=(B_t,\P_x)_{t\geq 0,x\in \R^d}$ 
be the Brownian motion in $\R^d (d\geq 1)$  independent of $S$. The subordinate Brownian motion is a stochastic process $X=(X_t,\P_x)_{t\geq 0,x\in \R^d}$ defined by $X_t=B_{S_t}$, $t\geq 0$\,. 

Recall that it is a L\' evy process with the characteristic (L\' evy) exponent $\phi(|\xi|^2)$ and it has transition density given by 
\begin{equation}\label{eq:heat_kernel}
p(t,x,y)=p(t,y-x)=\int_{(0,\infty)}(4\pi s)^{-d/2}e^{-\frac{|y-x|^2}{4s}}\P(S_t\in ds)
\end{equation}
for $x,y\in \R^d$ and $t>0$\,.

The following observation is important in obtaining estimates of transition density. Since $B$ and $S$ are independent, we may rewrite $\P_0(|X_t|\geq r)$, for $t,r>0$, as 
\begin{align*}
\P_0(|X_t|\geq r)&=\int_{(0,\infty)}\P_0(|B_s|^2\geq r^2)\P(S_t\in ds)\\
&=\int_{(0,\infty)}\P_0\left((2s)^{-1}|B_s|^2\geq (2s)^{-1}r^2\right)\P(S_t\in ds)\\
\end{align*}
Using the fact that $\frac{|B_s|^2}{2s}$ has chi-square distribution with parameter $d$ (as a sum of squares of $d$ independent standard normal random variables) and since it is independent of $S$, we conclude that 
\begin{align}
   \P_0(|X_t|\geq r)&=\int_{(0,\infty)}\P_0\left(Y\geq (2s)^{-1}r^2\right)\P(S_t\in ds)\nonumber\\
   &=\int_{(0,\infty)}\P(S_t\geq (2y)^{-1}r^2)\P_0(Y\in dy)\,,\label{eq:sbm-tail}
\end{align}
where
\begin{equation}\label{eq:Y}
\P_0(Y\geq t)=\frac{1}{2^{d/2}\Gamma(d/2)}\int_t^\infty s^{d/2-1}e^{-s/2}\,ds,\ t>0
\end{equation}
and $\Gamma(t):=\int_0^\infty y^{t-1}e^{-y}\,dy,\ t>0$ is the gamma function.

\subsection{Upper bounds}
We start with a lemma that shows that the inverse function of $\phi$ also satisfies certain scaling properties. 

\begin{Lem}\label{lem:inverse}
Let $\phi$ be a Bernstein function. 
\begin{itemize}
\item[(i)] Let $l:=\lim\limits_{\lambda\to 0+}\phi(\lambda)$ and $u:=\lim\limits_{\lambda\to\infty}\phi(\lambda)$. Then for any $l<\lambda<u$ and $x\geq 1$ such that $\lambda x<u$ we have $\frac{\phi^{-1}(\lambda x)}{\phi^{-1}(\lambda)}\geq x$\,.
\item[(ii)] If {\bf (L)} holds for $\phi$, then
\[
	\frac{\phi^{-1}(\lambda x)}{\phi^{-1}(\lambda)}\leq C_L^{-1/\gamma} x^{1/\gamma}\qquad \text{ for all }\quad \lambda>\phi(\lambda_L),\ x\geq C_L\,.
\]
\end{itemize}
\end{Lem}
\proof
(i) Since $\phi(\lambda x)\leq \phi(\lambda) x$ for any Bernstein function $\phi$ (see Lemma \ref{lem:bf}), $\lambda>0$ and $x\geq 1$, we can rewrite this (using the fact that $\phi$ is strictly increasing) as 
\[
    \frac{\phi^{-1}(\phi(\lambda)x)}{\phi^{-1}(\phi(\lambda))}\geq x\,.
\]
Setting $\eta=\phi(\lambda)$ we obtain $\frac{\phi^{-1}(\eta x)}{\phi^{-1}(\eta)}\geq x$ 
for all $l<\eta<u$ such that $\eta x<u$\,.\\
(ii)
Let $\lambda>\lambda_L$ and $x\geq 1$. Since $\phi^{-1}$ is increasing, one can rewrite {\bf (L)} as
\[
	\frac{\phi^{-1}(C_Lx^\gamma\phi(\lambda))}{\phi^{-1}(\phi(\lambda))}\leq x\,.
\]
Taking $\eta=\phi(\lambda)$ and $y=C_Lx^\gamma$ it follows that \[\frac{\phi^{-1}(\eta y)}{\phi^{-1}(\eta)}\leq C_L^{-1/\gamma} y^{1/\gamma},\ \eta>\phi(\lambda_L),\ y\geq C_L\,.\]
\qed

We proceed with the proof of the near diagonal upper bound for the transition density. This type of estimate can be found in \cite{KSch} for more general L\' evy processes via Nash inequality. Our estimate is obtained by using scaling properties of the L\' evy exponent. 
\begin{Prop}\label{prop:on-up}
	Assume that {\bf (L)} holds for $\phi$. There exists a constant $C_1>0$ such that 
	\[
		p(t,x)\leq C_1\left[\phi^{-1}(\tfrac{1}{t})\right]^{d/2}
	\]
	for all $0<t<\phi(\lambda_L)^{-1}$ and $x\in \R^d$\,. 
\end{Prop}
\proof
Since $e^{-t\phi(|\xi|^2)}=\int_{\R^d} e^{i\xi\cdot x}p(t,x)\,dx$,  the Fourier inversion formula yields
	\begin{equation}\label{eq:inv_fourier}
		p(t,x)=(2\pi)^{-d}\int_{\R^d}e^{-i\xi\cdot x}e^{-t\phi(|\xi|^2)}\,d\xi\leq \int_{\R^d} e^{-t\phi(|\xi|^2)}\,d\xi\,.
	\end{equation}
	
Switching to polar coordinates we get
\begin{align*}
	p(t,x)&\leq c_1\int_0^\infty r^{d-1}e^{-t\phi(r^2)}\,dr=c_1\int_0^\infty r^{d-1}\int_{t\phi(r^{2})}^\infty e^{-s}\,ds\,dr\\
	&=c_1\int_0^\infty e^{-s}\int_0^{\phi^{-1}(\frac{s}{t})^{1/2}}r^{d-1}\,dr\,ds= \tfrac{c_1}{d}\int_0^\infty e^{-s}\left[\phi^{-1}(\tfrac{s}{t})\right]^{\frac{d}{2}}\,ds
\end{align*}

Now we split the integral in the last display and use monotonicity of $\phi^{-1}$ and Lemma \ref{lem:inverse} again to obtain
\begin{align*}
	p(t,x)&\leq  \tfrac{c_1}{d}\left\{\int_0^1 e^{-s}\left[\phi^{-1}(\tfrac{s}{t})\right]^\frac{d}{2}\,ds+\int_1^\infty e^{-s}\left[\phi^{-1}(\tfrac{s}{t})\right]^\frac{d}{2}\,ds\right\}\\
	&\leq  \tfrac{c_1}{d}\left[\phi^{-1}(\tfrac{1}{t})\right]^\frac{d}{2}\left\{\int_0^1 e^{-s}\,ds+\int_1^\infty e^{-s}\left[\tfrac{\phi^{-1}(\frac{s}{t})}{\phi^{-1}(\frac{1}{t})}\right]^\frac{d}{2}\,ds\right\}\\
	&\leq  \tfrac{c_1}{d}\left[\phi^{-1}(\tfrac{1}{t})\right]^\frac{d}{2}\left\{\int_0^1 e^{-s}\,ds+C_L^{-\frac{d}{2\gamma}}\int_1^\infty e^{-s}s^{\frac{d}{2 \gamma }}\,ds\right\}\,.
\end{align*}
It is enough to take $C_1:= \tfrac{c_1}{d}\left\{\int_0^1 e^{-s}\,ds+C_L^{-\frac{d}{2\gamma}}\int_1^\infty e^{-s}s^{\frac{d}{2 \gamma }}\,ds\right\}$.
\qed

The next result is the off-diagonal upper bound. We start with an estimate of the tail of the chi-square distribution given by (\ref{eq:Y})\,.
\begin{Lem}\label{lem:chi-sq} For any $t_0>0$ there exists a constant $c=c(t_0)\geq 1$ such that 
\[
c^{-1}t^{d/2-1}e^{-t/2}\leq \P_0(Y\geq t)\leq ct^{d/2-1}e^{-t/2},\quad t\geq t_0\,.
\]
\end{Lem}
\proof
    Using integration by parts,
    \begin{equation}\label{eq:tmp3217}
        f(t):=\int_t^\infty s^{d/2-1}e^{-s/2} \,ds=2t^{d/2-1}e^{-t/2}+\left(d-2\right)\int_t^\infty s^{d/2-2}e^{-s/2} \,ds\,.
    \end{equation}
    Assume first that $d\geq 2$. Since
    $
        \int_t^\infty s^{d/2-2}e^{-s} \,ds\leq t^{-1}f(t)\,,
    $
    we get $f(t)(1-t^{-1}d)\leq 2t^{d/2-1}e^{-t/2}$, yielding
    \[
        f(t)\leq 4t^{d/2-1}e^{-t/2},\quad t\geq 2d\,.
    \]
    Therefore, in this case, from (\ref{eq:tmp3217}) we conclude that, for some constant $c_1=c_1(t_0)\geq 4$, 
    \[
        2t^{d/2-1}e^{-t/2}\leq f(t)\leq c_1 t^{d/2-1}e^{-t/2},\quad t\geq t_0\,.
    \]
    Let $d=1$. By (\ref{eq:tmp3217}) we get $f(t)\leq 2t^{d/2-1}e^{-t/2}$\,. On the other hand, 
    \[
        f(t)\geq \int_t^{2t}s^{-1/2}e^{-s/2}\,ds\geq 2^{1/2}t^{-1/2}(e^{-t/2}-e^{-t})\geq 2^{1/2}(1-e^{-t_0/2})t^{-1/2}e^{-t/2}
    \]
    for $t\geq t_0$\,.
\qed
\begin{Prop}\label{prop:off-up}
	Assume that {\bf (L)} holds for $\phi$. There exist constants $C_2>0$ and $a_U>0$ such that 
	\[
		p(t,x)\leq C_2 \left[t|x|^{-d}H(|x|^{-2})+\phi^{-1}(t^{-1})^{d/2}e^{-a_U|x|^2\phi^{-1}(t^{-1})}\right]\,.
	\]
	for all $0<t<(2e\phi(\lambda_L))^{-1}$ 
	and $x\in \R^d$ satisfying $t\phi(|x|^{-2})\leq 1$\,.
\end{Prop}
\proof
	 Let $r>0$, $0<t<(2e\phi(\lambda_L))^{-1}$ 
	 be such that $t\phi(r^{-2})\leq 1$. We  first note that 
	 \begin{equation}\label{eq:comp20134}
	 r^2\phi^{-1}(t^{-1})\geq 1\qquad \text{ and }\qquad 2e\leq \frac{\phi^{-1}(t^{-1})}{\phi^{-1}(\frac{1}{2e}t^{-1})}\leq \left(2e\right)^{1/\gamma}C_L^{-1/\gamma}\,,
	 \end{equation}
where the second inequality follows from Lemma \ref{lem:inverse}\,. Let $\wt{B}=(\wt{B}_t)_{t\geq 0}$ be $(d+2)$-dimensional Brownian motion and denote by $\wt{X}=(\wt{X}_t)_{t\geq 0}$ corresponding subordinate Brownian motion. Let $\wt{Y}$ be a random variable with the chi-square distribution with parameter $d+2$\,.
	  Using  (\ref{eq:sbm-tail}) and  Proposition \ref{prop:subord_upper} it follows that
	 \begin{align}
	 \P(|\wt{X}_t|\geq r)&\leq \int_{t\phi(2yr^{-2})\leq \frac{1}{2e}}4C_StH(2yr^{-2})\P_0(\wt{Y}\in dy)\nonumber\\
	 &+\int_{t\phi(2yr^{-2})\geq \frac{1}{2e}}\P_0(\wt{Y}\in dy)\nonumber\\
	 &\leq 4C_S tH(r^{-2})\int_{(0,\infty)}(1+(2y)^2)\P_0(\wt{Y}\in dy)\label{eq:tmp313}\\&+\P_0\left(\wt{Y}\geq \tfrac{1}{2}r^2\phi^{-1}(\tfrac{1}{2e}t^{-1})\right),\nonumber
	 \end{align}
since $\frac{H(2yr^{-2})}{H(r^{-2})}\leq 1\vee  (2y)^2$ by Lemma \ref{lem:bf}\,.

Using Lemma \ref{lem:chi-sq} we obtain
\begin{equation}\label{eq:tmp2367}
    \P(|\wt{X}_t|\geq r)\leq c_1 \left(tH(r^{-2})+r^{d}\phi^{-1}(\tfrac{1}{2e}t^{-1})^{d/2}e^{-\frac{1}{4}r^2\phi^{-1}(\frac{1}{2e}t^{-1})}\right)\,\,,
\end{equation}
since  $\tfrac{1}{2}r^2\phi^{-1}(\tfrac{1}{2e}t^{-1})\geq \tfrac{1}{4}C_L^{1/\gamma}\left(2e\right)^{-1/\gamma}$ by (\ref{eq:comp20134})\,.

On the other hand, from (\ref{eq:heat_kernel}),
\begin{align*}
\P(|\wt{X}_t|\geq r)&=c_2\int_{|x|\geq r}\int_{(0,\infty)} s^{-(d+2)/2}e^{-\frac{|x|^2}{4s}}\,\P(S_t\in ds)\,dx\\
&=c_3\int_{(0,\infty)}s^{-d/2-1}\int_r^\infty y^{d+1}e^{-\frac{y^2}{4s}}\,dy\,\P(S_t\in ds)\,.
\end{align*}
Integration by parts yields
\begin{align*}
     \int_r^\infty y^{d+1}e^{-\frac{y^2}{4s}}\,dy&=2sr^{d}e^{-\frac{r^2}{4s}}+2ds\int_r^\infty y^{d-1}e^{-\frac{y^2}{4s}}\,dy\\
    &\geq 2sr^{d}e^{-\frac{r^2}{4s}}\,.
\end{align*}
Hence,
\[
      \P(|\wt{X}_t|\geq r)\geq 2c_3r^d\int_{(0,\infty)} s^{-d/2}e^{-\frac{r^2}{4s}}\P(S_t\in ds)=2c_3|x|^dp(t,x)\,,
\]
for $x\in \R^d$ such that $|x|=r$ by (\ref{eq:heat_kernel})\,. Using the  last display (\ref{eq:tmp2367}) and (\ref{eq:comp20134}) we get, for some $a_U>0$, 
\[
    |x|^dp(t,x)\leq c_4\left(tH(|x|^{-2})+|x|^{d}\phi^{-1}(t^{-1})^{d/2}e^{-a_U|x|^2\phi^{-1}(t^{-1})}\right)
\]
$0<t<(2e\phi(\lambda_L))^{-1}$ 
and $x\in \R^d$ satisfying $t\phi(|x|^{-2})\leq 1$\,.

\qed


\subsection{Lower bounds} 
In this subsection the lower bounds of the transition density will be proved.

\begin{Prop}\label{prop:on-low}
	Assume that {\bf (L)} holds for $\phi$. There exist constants $C_4>0$ and $\kappa\in (0,1)$  such that
\[
		p(t,x)\geq C_4\phi^{-1}(t^{-1})^{\frac{d}{2}}e^{-2|x|^{2}\phi^{-1}(t^{-1})}
	\]
	for all $0<t<\kappa\phi(\lambda_L)^{-1}$ and $x\in \R^d$\,.
	In particular, if additionally $t\phi(|x|^{-2})\geq 1$ holds, we have 
	\[
	    p(t,x)\geq C_5\phi^{-1}(t^{-1})^{\frac{d}{2}}\,,
	\]
	where $C_5=C_4e^{-2}$\,.
\end{Prop}
\proof First we note that {\bf (L)} implies $\lim\limits_{\lambda \to\infty}\phi(\lambda)=\infty$. 
Let $0<t<\phi(\lambda_L)^{-1}$ and  $x\in \R^d$
. Using (\ref{eq:heat_kernel}) we get
\begin{align}
p(t,x)&=c_1\int_{(0,\infty)}s^{-d/2}e^{-\frac{|x|^2}{4s}}\P(S_t\in ds)\nonumber \\
    &\geq c_1 \int_{\left[2^{-1}\phi^{-1}(t^{-1})^{-1},\phi^{-1}(ct^{-1})^{-1}\right]}s^{-d/2}e^{-\frac{|x|^2}{4s}}\P(S_t\in ds)\nonumber \\
  &\geq c_1 \phi^{-1}(\rho t^{-1})^{d/2}e^{-\frac12 |x|^2\phi^{-1}(t^{-1})}\P\left(2^{-1}\phi^{-1}(t^{-1})^{-1}\leq S_t\leq\phi^{-1}(\rho t^{-1})^{-1}\right),
\label{eq:tmp_2034}
\end{align}
where $\alpha=2, \beta=1$ and $\rho\in (0,1)$ are the constants from Proposition \ref{prop:sub-low}. 
By Lemma \ref{lem:inverse},  for $0<t<\rho\phi(\lambda_L)^{-1}$,
\[
     \phi^{-1}(\rho t^{-1})=\phi^{-1}(t^{-1})\frac{\phi^{-1}\rho t^{-1})}{\phi^{-1}(t^{-1})}\geq C_L^{1/\gamma}\rho^{1/\gamma} \,\phi^{-1}(t^{-1}).
\]
Using the last display and Proposition \ref{prop:sub-low} in (\ref{eq:tmp_2034}) we get
\[
    p(t,x)\geq c_1C_L^{d/2\gamma}\tau e^{-\frac12 |x|^{-2}\phi^{-1}(t^{-1})}\phi^{-1}(t^{-1})^{d/2}\,. 
\]
If $t\phi(|x|^{-2})\geq 1$, then $|x|^{2}\phi^{-1}(t^{-1})\leq 1$. 
Hence, from the last display we conclude
\[
     p(t,x)\geq c_1C_L^{d/2\gamma}\rho^{1/\gamma}\tau e^{-2}\phi^{-1}(t^{-1})^{d/2}\,. 
\]
\qed


\begin{Prop}\label{prop:off-lower2} 
Assume that {\bf (L)} and {\bf (U)} hold for $H$ with $\delta<2$ and that the drift of the subordinator is zero. There exist constants 
$C_6>0$, $L>1$ and $M\in (0,1)$  
such that for all $0<t<\frac{1}{2}\phi(\lambda_L)^{-1}$ and $x\in \R^d$ satisfying 
$|x|<\sqrt{L^{-1}\lambda_L^{-1}\wedge M\lambda_U^{-1}}$  
and $t\phi(|x|^{-2})\leq 1$ we have
 \[
  p(t,x)\geq C_6 t|x|^{-d}H(|x|^{-2})\,.
  \]
\end{Prop}
\proof 
From Proposition \ref{prop:subord-lower} it follows that there exist $M\in (0,1)$, $L>1$ and $c_S>0$ such that for $|x|^2<L^{-1}\lambda_L^{-1} \wedge M\lambda_U^{-1}$ it holds that 
\[
    \P(|x|^2\leq S_t\leq L |x|^2)\geq c_StH(|x|^{-2})\,.
\]
Hence, by (\ref{eq:heat_kernel}) we obtain
\begin{align*}
 p(t,x)&\geq c_1 \int_{[|x|^2,L|x|^2]}s^{-d/2}e^{-\frac{|x|^2}{4s}}\P(S_t\in ds)\\
 &\geq c_1 L^{-d/2}|x|^{-d}e^{-1/4}\P(|x|^2\leq S_t\leq L |x|^2)\\
 &\geq c_2 t|x|^{-d}H(|x|^{-2})\,.
 \end{align*}

\qed

\subsection{Proofs of main results}
In this section we prove our main results. We start with the proof of the equivalence of the lower scaling condition {\bf (L)} for $\phi$ and the near diagonal upper bound of the heat kernel. 

\proof[Proof of Theorem \ref{tm:upper-equiv}]
If {\bf (L)} holds for $\phi$, then we use Propositon \ref{prop:on-up}. Assume that (\ref{eq:on-up}) holds. Then 
\[
    \int_{(0,\infty)}(4\pi s)^{-d/2}\P(S_t\in ds)
    =p(t,0)\leq C_1\phi^{-1}(t^{-1})^{d/2}
\]
for $0<t<\phi(\lambda_L)^{-1}$. In particular, for $a>1$ we have
\begin{align*}
C_1\phi^{-1}(t^{-1})^{d/2}&\geq \int_{(0,2a\phi^{-1}(at^{-1})^{-1})}(4\pi s)^{-d/2}\P(S_{t}\in ds)\\&\geq (8\pi a)^{-d/2}\phi^{-1}(at^{-1})^{d/2}\P(S_t< 2a\phi^{-1}(at^{-1})^{-1})\,.
\end{align*}
Since
\begin{align*}
    \P(S_{t}< 2a\phi^{-1}(at^{-1})^{-1})&=1-\P(1-e^{-\phi^{-1}(at^{-1})S_{t}}\geq 1-e^{-2a})\\&\geq 1-\frac{1-e^{-t\phi(\phi^{-1}(at^{-1}))}}{1-e^{-2a}}=\frac{e^{-a}-e^{-2a}}{1-e^{-2a}}=:c_1>0
\end{align*}
we obtain
\[
    \frac{\phi^{-1}(at^{-1})}{\phi^{-1}(t^{-1})}\leq 8\pi a(C/c_1)^{2/d}=:b,\quad 0<t<\phi(\lambda_L)^{-1}\,.
\]
Note that $b>1$. 
Taking $\lambda=\phi^{-1}(t^{-1})$ in the last display we get
\begin{equation}\label{eq:tmp-upper-1step}
    \frac{\phi(\lambda b)}{\phi(\lambda)}\geq a,\qquad \lambda>\lambda_L\,.
\end{equation}  

Let $x\geq 1$ and choose $n\in \N$ so that $b^{n-1}\leq x<b^n$. By iterating (\ref{eq:tmp-upper-1step}), for any $\lambda>\lambda_L$ we get
\[
    \frac{\phi(\lambda x)}{\phi(\lambda)}\geq \frac{\phi(\lambda b^{n-1} )}{\phi(\lambda)}\geq a^{n-1}=b^{(n-1)\frac{\log{a}}{\log{b}}}\geq b^{-\frac{\log a}{\log b}}x^{\frac{\log a}{\log b}}\,.
\]
Hence, {\bf (L)} holds with $\gamma=\frac{\log a}{\log b}$. 
\qed

Proof of the main result is a consequence of propositions established in the previous subsections. 
\proof[Proof of Theorem \ref{tm:main}]
  Follows directly from Proposition \ref{prop:on-up}, Proposition \ref{prop:on-low}, Proposition \ref{prop:off-up}  and Proposition \ref {prop:off-lower2}  since {\bf(L)} also holds for $\phi$ by Lemma \ref{lem:bf}. 
\qed

\proof[Proof of Corollary \ref{cor:classical}]
Note first that, by Proposition \ref{prop:comp}, $\phi$ and $H$ are comparable, so {\bf (L)} (resp. {\bf (U)}) are basically the same for these two functions. 

Assume that $t\phi(|x|^{-2})\leq 1$.  
Then $|x|^2\phi^{-1}(t^{-1})\geq 1$
and so there exists a constant $c_1\geq 1$ so that 
\begin{equation*}\label{eq:expineq}
 e^{-a_U|x|^2\phi^{-1}(t^{-1})}\leq c_1\left(|x|^2\phi^{-1}(t^{-1})\right)^{-d/2-1}\,.
\end{equation*}
Hence, by Lemma \ref{lem:inverse} (i), 
\begin{align}\label{eq:comp9999}
\phi^{-1}(t^{-1})^{d/2}e^{-a_U|x|^2\phi^{-1}(t^{-1})}&\leq c_1|x|^{-d}\frac{\phi^{-1}(\phi(|x|^{-2}))}{\phi^{-1}\left(\frac{\phi(|x|^{-2})}{t\phi(|x|^{-2})}\right)}\leq c_1 |x|^{-d}t\phi(|x|^{-2})\,,
\end{align}
implying the upper bound 
\begin{align*}
  t|x|^{-d}H(|x|^{-2})\vee \phi^{-1}(t^{-1})^{d/2}e^{-a_U|x|^2\phi^{-1}(t^{-1})}&\leq t|x|^{-d}\phi(|x|^{-2})\vee c_1t|x|^{-d}\phi(|x|^{-2})\\&\leq c_1t|x|^{-d}\phi(|x|^{-2})\,.
\end{align*}
On the other hand, by Proposition \ref{prop:comp},
\[
 t|x|^{-d}H(|x|^{-2})\vee \phi^{-1}(t^{-1})^{d/2}e^{-a_L|x|^2\phi^{-1}(t^{-1})}\geq t|x|^{-d}H(|x|^{-2})\geq c_2 t|x|^{-d}\phi(|x|^{-2})
\]
for $|x|^{-2}>\lambda_U$\,.
The estimate follows now from Theorem \ref{tm:main}\,.
\qed

It is left to prove the Green function estimates using the heat kernel estimates. 

\proof[Proof of Corollary \ref{cor:green}]
 We split the integral
 \[
  G(x)=\int_0^{\phi(|x|^{-2})^{-1}}p(t,x)\,dt+\int_{\phi(|x|^{-2})^{-1}}^\infty p(t,x)\,dt
 \]
 and use already established upper and lower bounds. 
 By Proposition \ref{prop:on-low} and  Lemma \ref{lem:inverse},
 \begin{align*}
 \allowdisplaybreaks
  \int_{\phi(|x|^{-2})^{-1}}^\infty p(t,x)\,dt&\geq c_1\int_{\phi(|x|^{-2})^{-1}}^{2\phi(|x|^{-2})^{-1}} \phi^{-1}(t^{-1})^{d/2}\,dt\\
  &\geq c_1 \phi^{-1}\left(\frac{\phi(|x|^{-2})}{2}\right)^{d/2}\cdot \frac{1}{\phi(|x|^{-2})}
  \geq \frac{c_1C_L^{1/\gamma}2^{-1/\gamma}}{|x|^d\phi(|x|^{-2})}\,.
 \end{align*}
 On the other hand, using the near diagonal upper bound from Proposition \ref{prop:on-up}  we get
  \begin{align*}
 \allowdisplaybreaks
  \int_{\phi(|x|^{-2})^{-1}}^\infty p(t,x)\,dt&\leq c_2\int_{\phi(|x|^{-2})^{-1}}^\infty \phi^{-1}(t^{-1})^{d/2}\,dt=c_2\int_0^{|x|^{-2}}y^{d/2}\left(-\frac{1}{\phi(y)}\right)'\,dy\\
  &=c_2\int_0^{|x|^{-2}}y^{d/2-1}\frac{y\phi'(y)}{\phi(y)^2}\,dy\leq c_2\int_0^{|x|^{-2}}\frac{y^{d/2-1}}{\phi(y)}\,dy\,,
  \end{align*}
  since $\lambda\phi'(\lambda)\leq \phi(\lambda)$\,. If $d\geq 3$ we take $\delta=1$ by Lemma \ref{lem:bf} and for $d=1,2$ we use {\bf (U)} and the assumption $2\delta<d$. In both cases we have $\frac{d}{2}-1-\delta>-1$\,. Hence
  \[
      \int_0^{|x|^{-2}}\frac{y^{d/2-1}}{\phi(y)}\,dy\leq \frac{c_3}{\phi(|x|^{-2})}\int_0^{|x|^{-2}} y^{d/2-1}\left(\frac{|x|^{-2}}{y}\right)^{\delta}\,dy=\frac{c_4}{|x|^d\phi(|x|^{-2})}
  \]
  
 It is left to estimate the integral of the off-diagonal upper bound from Proposition \ref{prop:off-up} integrated over the appropriate time-interval:
 \begin{align*}
  \int_0^{\phi(|x|^{-2})^{-1}} p(t,x)\,dt&\leq c_4 \int_0^{\phi(|x|^{-2})^{-1}}\left(t|x|^{-d}H(|x|^{-2})+\phi^{-1}(t^{-1})^{d/2}e^{-a_U|x|^2\phi^{-1}(t^{-1})}\right)\,dt\\
  &\leq c_5 (|x|^{-d}H(|x|^{-2})+c_4|x|^{-d}\phi(|x|^{-2}))\int_0^{\phi(|x|^{-2})^{-1}}t\,dt \\
  &=\tfrac{1}{2}c_5 |x|^{-d}\left(\frac{H(|x|^{-2})}{\phi(|x|^{-2})^2}+\frac{1}{\phi(|x|^{-2})}\right)\leq \frac{c_5}{|x|^d\phi(|x|^{-2})}\,,
 \end{align*}
 where in the second inequality we have used (\ref{eq:comp9999}) and in the last $H(\lambda)\leq \phi(\lambda)$\,.

\qed
\providecommand{\bysame}{\leavevmode\hbox to3em{\hrulefill}\thinspace}
\providecommand{\MR}{\relax\ifhmode\unskip\space\fi MR }
\providecommand{\MRhref}[2]{%
  \href{http://www.ams.org/mathscinet-getitem?mr=#1}{#2}
}
\providecommand{\href}[2]{#2}

\end{document}